\documentstyle[pb-diagram]{article}

\newcommand{\es}{\emptyset}

\newcommand{\ba}{\begin{array}}
\newcommand{\ea}{\end{array}}

\newtheorem{theorem}{Theorem}
\newtheorem{proposition}[theorem]{Proposition}
\newtheorem{lemma}[theorem]{Lemma}

\newtheorem{definition}[theorem]{Definition}
\newtheorem{corollary}[theorem]{Corollary}

\newcommand{\be}{\begin{enumerate}}
\newcommand{\ee}{\end{enumerate}}
\newcommand{\bi}{\begin{itemize}}
\newcommand{\ei}{\end{itemize}}
\newcommand{\bd}{\begin{description}}
\newcommand{\ed}{\end{description}}

\newcommand{\et}{\wedge}

\newcommand{\vel}{\vee}

\newcommand{\imp}{\rightarrow}

\newcommand{\beq}{\begin{eqnarray*}}
\newcommand{\eeq}{\end{eqnarray*}}
\newcommand{\seq}{\Rightarrow}
\author{ {F.Parlamento, F.Previale }
\\Department of Mathematics,  Computer Science and Physics
\\University of Udine,  via  Delle Scienze 206, 33100 Udine, Italy.
\\Department of Mathematics
\\University of Turin, via Carlo Alberto 10, 10123 Torino, Italy
\\e-mail: {\em franco.parlamento$@$uniud.it} 
}
 \title{Absorbing the Structural Rules in the Sequent Calculus with Additional Atomic Rules}
\date{}

\begin{document}
\maketitle
%\newpage

%\tableofcontents

%\newpage

%\noindent{\bf Mathematical Subject Classification: 03F05}

\begin{abstract}
We  show that if the structural rules are admissible over a set ${\cal R}$ of atomic rules, then they are admissible in the sequent calculus obtained by adding the rules in ${\cal R}$ to ${\bf G3[mic]}$.
Two applications to pure logic and to the sequent calculus with equality are  presented.

\end{abstract}


\


\section{Introduction}

\label{intro}

A multisuccedent sequent calculus  for intuitionistic logic free of structural rules was presented   in \cite{D88} and a
detailed proof of their admissibility,   based on \cite{D97},  appeared in \cite{NvP01}. A single succedent version of that calculus
was adopted in \cite{TS96}. In all cases the proof of the admissibility of the structural rules relies, as for the classical $ {\bf G3}$ system, on the hight-preserving admissibility of the contraction rule.
When additional atomic rules are added to the calculus the hight-preserving admissibilibility of the contraction rule may fail. Such is the case for example for the following rules 
 $\hbox{Ref}$ and $ \hbox{Repl}$   for equality, introduced in \cite{NvP98} and adopted in  the second edition \cite{TS00} of \cite{TS96}
    \[
\ba{clccl}
\underline{t=t,\Gamma\seq \Delta}&\vbox to 0pt{\hbox{~Ref}}&~~~~~~~~~~~~~~~~~~~~~&\underline{s=r, P[x/s],P[x/r], \Gamma \seq \Delta}&\vbox to 0pt{\hbox{~Repl}}\\
\Gamma\seq \Delta&&&s=r, P[x/s], \Gamma \seq \Delta&
\ea
\]

For example 
\[
a=f(a), a=f(a) \seq  a=f(f(a))
\]
has  derivations of height equal 1  in the systems obtained by adding $\hbox{Ref}$ and $\hbox{Repl}$  to ${\bf G3[mic]}$ 
 namely 
\[
\ba{c}
a=f(a), a=f(a),  a=f(f(a)) \seq a=f(f(a))\\
\cline{1-1}
a=f(a), a=f(a) \seq  a=f(f(a))
\ea
\]
but $a=f(a)\seq a=f(f(a))$ cannot have  a derivation of height less than or equal 1 in such a system.

For that reason, in such cases, to prove the admissibility of the structural rules    we  have to  follow a route somewhat different from
the one used in  \cite{NvP98} (see also \cite{NvP01} and \cite{NvP11}),
and from the one followed in \cite{TS00} for extensions of their single succedent ${\bf G3[mic]}$ calculus with  rules for which hight preserving admissibilty of the contraction rule is ensured. The basic idea is to eliminate context-sharing cuts first,  with the eliminability of contraction obtained as a consequence, due to its immediate derivability from context-sharing cut.

Actually we  will show that we can proceed in that way for any set ${\cal R}$ of atomic formulae of the following form:

\[
\ba{cl}
\underline{\vec{Q_1},\Gamma_1\seq \Delta_1, \vec{Q_1'} ~~\ldots ~~\vec{Q_n},\Gamma_n\seq \Delta_n, \vec{Q_n'}}&\\
\vec{P}, \Gamma_1,\ldots, \Gamma_n \seq \Delta_1,\ldots, \Delta_n, \vec{P'}
\ea
\]
where $\vec{Q_1}, \vec{Q_1'} ,\ldots, \vec{Q_n},  \vec{Q_n'}, \vec{P},\vec{P'}$ are 
 sequences (possibly empty) of  atomic formulae and $\Gamma_1,\ldots, \Gamma_n, \Delta_1,\ldots \Delta_n$  are  finite sequences (possibly empty) of formulae that are not active in the rule.
 
More precisely,  letting  ${\bf G3[mic]}^{\cal R}$ denote the calculi  obtained by adding to ${\bf G3[mic]}$ the rules in ${\cal R}$  and $Cut_{cs}$   the context-sharing cut rule,  we will show that
 any derivation in  ${\bf G3[mic]}^{\cal R} + Cut_{cs}$ can be transformed into a derivation in the same system
 in which the  rules in ${\cal R}$ and the $Cut_{cs}$ rule are applied before any logical rule. From that it will follow that
if the structural rules are admissible  in the calculus that contains only the initial sequents  and the rules in ${\cal R}$, then they are admissible in  ${\bf G3[mic]}^{\cal R}$ as well.
For ${\cal R}=\es$  we have   that  the height preserving admissibility of the weakening rules and the height-preserving invertibility of the logical rules  suffice for the eliminability of context-sharing cut, without having to obtain the admissibility of contraction first, and, as a  consequence, for the admissibility of the $Cut$-rule in ${\bf G3[mic]}$.
For ${\cal R}= \{ \hbox{Ref},   \hbox{Repl}\}$ we obtain  that the structural rules are admissible in ${\bf G3[mic]}^{\cal R}$, thus extending the result proved in \cite{NvP98} in the case $t,r$ and $s$ are restricted to be constants.

\section{Preliminaries }
\label{Preliminaries}

The sequent calculus ${\bf G3c}$ in \cite{N16} has, in  the notations in \cite{TS00},  the following initial sequents and rules, where $P$ is an atomic formula and $A, B$ stand for any formula in a first order language   (function symbols included) and $\Gamma$ and $\Delta$ are finite multisets of formulae :

\

{\bf Initial sequents}
\[
P,\Gamma \seq \Delta, P
\]

{\bf Logical rules}

\[
\ba{clccl}
A,B, \Gamma \seq \Delta&\vbox to 0pt{\hbox{$L\et $}}&~~~~~~&\underline{\Gamma\seq,\Delta, A~~~\Gamma\seq \Delta, B}& \vbox to 0pt{\hbox{$R\et$}}\\
%\cline{1-1}\cline{4-4}
\overline{A\et B, \Gamma \seq \Delta}&&&\Gamma \seq \Delta, A\et B\\

\\
\underline{A,\Gamma\seq \Delta~~~B,\Gamma \seq \Delta}& \vbox to 0pt{\hbox{$L\vel $}}&&\Gamma\seq \Delta, A, B& \vbox to 0pt{\hbox{$R\vel$}}\\
%\cline{1-1}\cline{4-4}
A\vel B, \Gamma \seq \Delta&&&\overline{\Gamma\seq\Delta, A\vel B}\\

\\

\\

\underline{\Gamma\seq \Delta, A~~~B,\Gamma \seq \Delta}&\vbox to 0pt{\hbox{$L\imp $}}&&A,\Gamma \seq \Delta, B&\vbox to 0pt{\hbox{$R\imp $}}\\
%\cline{1-1}\cline{4-4}
A\imp B, \Gamma \seq \Delta&&&\overline{\Gamma \seq \Delta, A\imp B}&\\
\\

\\
\overline{\bot, \Gamma \seq \Delta}^{~L\bot}&&&&\\

\\
\\

\underline{A[x/t], \forall x A,\Gamma\seq \Delta}&\vbox to 0pt{\hbox{$L\forall $}}&&\underline{\Gamma\seq \Delta, A[x/a]}&\vbox to 0pt{\hbox{$R\forall $}}\\
%\cline{1-1}\cline{4-4}
\forall xA, \Gamma\seq \Delta&&&\Gamma\seq \Delta, \forall x A\\

\\

\\
\underline{A[x/a], \Gamma\seq \Delta}&\vbox to 0pt{\hbox{$L\exists $}}&&\underline{\Gamma\seq \Delta, \exists x A,  A[x/t]}&\vbox to 0pt{\hbox{$R\exists $}}\\
%\cline{1-1}\cline{4-4}
\exists xA, \Gamma\seq \Delta&&&\Gamma\seq \Delta, \exists x A

\ea
\]

In ${\bf G3i}$ the rules $L\imp$, $R\imp$ and $R\forall$ are replaced by:

\[
\ba{clccl}
\underline{A\imp B, \Gamma \seq \Delta, A~~~B,\Gamma \seq \Delta}&\vbox to 0pt{\hbox{$L^i\imp $}}&~~~~~~~~&A,\Gamma \seq B&\vbox to 0pt{\hbox{$R^i\imp $}}\\
%\cline{1-1}\cline{4-4}
A\imp B, \Gamma\seq \Delta&&&\overline{\Gamma\seq\Delta, A\imp B}\\

\\
\\

&&&\Gamma\seq A[x/a]&\vbox to 0pt{\hbox{$R^i\forall $}}\\
%\cline{4-4}
&&&\overline{\Gamma\seq \Delta, \forall x A}&

\ea
\]
In  both ${\bf G3c}$ and ${\bf G3i}$,   $a$ does not occur in the conclusion of $L\exists$ and $R\forall$.

Finally ${\bf G3m}$ is obtained from ${\bf G3i}$ by replacing $L\bot$ by the initial sequents $\bot, \Gamma \seq \Delta, \bot$.

\

${\bf G3[mic]}$ denotes any of the systems ${\bf G3m}$,   ${\bf G3i}$ or  ${\bf G3c}$.

\

The left and right weakening rules, $LW$ and $RW$ have the form:

\[
\ba{clccl}
\Gamma\seq \Delta&\vbox to 0pt{\hbox{$LW$}}&~~~~&\Gamma\seq \Delta&\vbox to 0pt{\hbox{$RW$}}\\
\cline{1-1}\cline{4-4}
A,\Gamma\seq \Delta&&&\Gamma\seq\Delta,A
\ea
\]
The left and right contraction rules, $LC$ and $RC$ have the form:

\[
\ba{clccl}
A.A, \Gamma\seq \Delta&\vbox to 0pt{\hbox{$LC$}}&~~~~&\Gamma\seq \Delta, A,A&\vbox to 0pt{\hbox{$RC$}}\\
\cline{1-1}\cline{4-4}
 A,\Gamma\seq \Delta&&&\Gamma\seq\Delta,A
 \ea
\]

The cut rule  and the context-sharing cut rule, $Cut$  and $Cut_{cs}$  have the form:
\[
\ba{clcclccl}
\Gamma\seq\Delta, A~~~~A,\Lambda \seq \Theta&\vbox to 0pt{\hbox{$Cut $}}&~~&&~~&&\Gamma\seq \Delta, A~~~A,\Gamma\seq \Delta&\vbox to 0pt{\hbox{$Cut_{cs} $}}\\
\cline{1-1}\cline{7-7}
\Gamma,\Delta\seq \Lambda,\Theta&&&&&&\Gamma\seq \Delta\\
\ea
\]

\

The additional atomic rules that we will consider are  of the form described in  \ref{intro}.

\subsection{Separated derivations}
\label{sep}

In all the systems considered the weakening rules are height-preserving admissible. For left weakening   it suffices to add $A$ to the antecedent  of every sequent in a given derivation of $\Gamma \seq \Delta$, modulo a possible renaming of the proper variables in the $L\exists$ and $R\forall$ inferences,  to obtain a derivation of the same height of $A, \Gamma\seq \Delta$. In the classical case one can proceed in the same way also for right weakening, while in the minimal or intuitionistic case one uses  induction on the height of derivation, taking advantage of the possibility of adding an arbitrary context on the right in the applications of $R^i\imp$ and $R^i\forall$.

\begin{definition}
For a set of additional  atomic rules ${\cal R}$, let ${\bf G3[mic]}^{\cal R}$ be the sequent calculus obtained from  ${\bf G3[mic]}$
by adding  the rules in ${\cal R}$. 
With ${\cal R}$ we will denote also the  logic-free subcalculus of  ${\bf G3[mic]}^{\cal R}$, that contains only the initial sequents $P,\Gamma\seq \Delta, P$  and the rules in ${\cal R}$.
\end{definition}

\begin{definition}
Let ${\cal R}$ be any set of atomic rules. An ${\cal R}$-inference is any application of a rule in ${\cal R}$. An ${\cal R}$-derivation is a derivation in 
${\bf G3[mic]}^{\cal R}$.
\end{definition}

\begin{definition}
A derivation in  ${\bf G3[mic]}^{\cal R} +Cut_{cs}$ 
 is  said to be {\em separated} if no logical inference  precedes an ${\cal R}$  or $Cut_{cs}$-inference.
\end{definition}

Our first goal  is to show that every derivable sequent in ${\bf G3[mic]}^{\cal R} +Cut_{cs}$ 
has a separated derivation in the same system. 
Derivations without  logical inferences are trivially separated. For such derivations we have the following useful fact.

\begin{lemma}\label{atomic}
If $\Gamma\seq \Delta$ has a derivation ${\cal D}$ in ${\cal R}$, then there is a subsequent $\Gamma^\circ\seq \Delta^\circ$ of $\Gamma\seq \Delta$, with atomic formulae only, that has a derivation ${\cal D}^\circ$ in the same system, with only atomic sequents, such that $h({\cal D}^\circ)\leq h({\cal D})$.
\end{lemma}

{\bf Proof} The claim is  proved by a straightforward induction on the height of derivations, thanks to the height-preserving admissibility of the weakening rules. $\Box$

\begin{lemma}\label{invertzero}
If the conclusion of a classical logical rule has a derivation in ${\cal R}$  of height bounded by $h$, then also its premisses have  derivations
 in ${\cal R}$ of height bounded by $h$. The same holds in the minimal  and intuitionistic case with the exception of rules $R\imp$ and $R\forall$.
\end{lemma}
{\bf Proof} By the previous lemma, there is a derivation  ${\cal D}^\circ$ of height bounded by $h$ of an atomic subsequent $\Gamma^\circ\seq \Delta^\circ$ of the conclusion of the logical inference. Being  atomic  $\Gamma^\circ\seq \Delta^\circ$  does not contain the principal formula of the logical inference we are interested in and its premiss or premisses can be obtained from  ${\cal D}^\circ$   by weakening. $\Box$

\begin{proposition}\label{invert}
\bi
\item[a)]  Hight-preserving separated invertibility of the logical rules in ${\bf G3c}^{\cal R} +Cut_{cs}$ 

If the conclusion of a logical inference has a separated derivation of height bounded by $h$, then also its  premisses  have  separated derivations of height bounded by $h$.

\item[b)] The same holds for ${\bf G3[mi]}^{\cal R} +Cut_{cs}$, except for the rules  $R^i\imp$ and $R^i\forall$.
\ei

\end{proposition}

{\bf Proof}  If the given derivation ${\cal D}$ reduces to an initial sequent or to an instance of $L\bot$ or it ends with a logical inference that does not introduce the principal formula of the rule $R$ to be proved invertible, then the argument is the same as for  the ${\bf G3[mic]}$ systems (see \cite{D88},  \cite{NvP01}).
For example if ${\cal D}$ has the form:

\[
\ba{ccc}
{\cal D}_0&~~~&{\cal D}_1\\
E\imp F, A\imp B, \Gamma \seq \Delta, E&&F, A\imp B, \Gamma\seq \Delta\\
\cline{1-3}
\multicolumn{3}{c}{A\imp B, E\imp F, \Gamma\seq \Delta}
\ea
\]
and the rule to be proved invertible is an $L^i\imp$ with principal formula $A\imp B$, then a derivation of the same height as ${\cal D}$ of
its first premiss $A\imp B, E\imp F, \Gamma\seq \Delta, A$ is obtained by height preserving weakening applied to ${\cal D}$.
As far as the second premiss, i.e. $B, E\imp F, \Gamma \seq \Delta$ is concerned, by  induction hypothesis there is a separated derivation ${\cal D}_0'$ of height bounded by the height of ${\cal D}_0$, of $B, E\imp F, \Gamma \seq \Delta, E$ and a separated derivation ${\cal D}_1'$,  of height bounded by the height of ${\cal D}_1$, of $B, F,  \Gamma \seq \Delta$, Then:
\[
\ba{ccc}
{\cal D}_0'&~~~&{\cal D}_1'\\
B, E\imp F, \Gamma\seq \Delta, E&&B, F,  \Gamma\seq \Delta\\
\cline{1-3}
\multicolumn{3}{c}{B,  E\imp F, \Gamma\seq \Delta}

\ea
\]is a separated derivation of   $B, E\imp F, \Gamma \seq \Delta$ with height  bounded by $h({\cal D})$.

 If the last inference does introduce  the principal formula of $R$, we only need, in addition, to note that the subderivations of a separated derivation are themselves  separated. 
 
 If ${\cal D}$ ends with an ${\cal R}$ or  $Cut_{cs}$-inference,  then ${\cal D}$, being separated, does not contain any logical inference, and the previous Lemma applies. $\Box$

\begin{proposition} \label{sepequality}
If the premisses  of an ${\cal R}$-inference $R$ have  a separated derivation, then its conclusion also has a separated derivation.

\end{proposition}

{\bf Proof} If all the separated derivations of the premisses end with an ${\cal R}$ or with a $Cut_{cs}$-inferences, they are all free of logical inferences and it suffices to apply to such premisses the rule $R$ to obtain the desired separated derivation.
Otherwise we select a derivation of a premiss that ends with a logical inference and proceed by a straightforward induction on the sum of the heights of derivations. For example suppose  that $R$ has the following two premisses:
1) $Q_1,\Gamma_1', A\vel B\seq \Delta_1, Q_1'$, with separated derivation ${\cal E}_1$, that ends with an $L\vel$-inference with principal formula $A\vel B$, and  2) $Q_2,\Gamma_2\seq \Delta_2, Q_2'$ with separated derivation ${\cal E}_2$
and that the  conclusion of $R$ is  $P,\Gamma_1',A\vel B, \Gamma_2\seq \Delta_1,\Delta_2, P'$.
${\cal E}_1$ has the form:

\[
\ba{ccc}
{\cal E}_{10}&~~~&{\cal E}_{11}\\
Q_1, \Gamma_1', A\seq \Delta_1, Q_1'&&Q_1, \Gamma_1', B\seq \Delta_1, Q_1'\\
\cline{1-3}
\multicolumn{3}{c}{Q_1,\Gamma_1', A\vel B \seq \Delta_1, Q_1'}
\ea
\]

We can apply the induction hypothesis to ${\cal E}_{10}$ paired with ${\cal E}_2$ to obtain a separated derivation of:
$$a)~~~P,\Gamma_1', A, \Gamma_2\seq \Delta_1,\Delta_2, P'$$ and  to ${\cal E}_{11}$ paired with ${\cal E}_2$ to obtain a separated derivation of:
$$b)~~P,\Gamma_1', B, \Gamma_2\seq \Delta_1,\Delta_2, P'.$$
Then it suffices to apply the same last  $L\vel$-inference of ${\cal E}_1$ to $a)$ and $b)$ to obtain the desired separated derivation of 
$$P,\Gamma_1',A\vel B, \Gamma_2\seq \Delta_1,\Delta_2, P'.$$
$\Box$

\begin{lemma}
If the premisses $\Gamma \seq \Delta, A$  and $A, \Gamma\seq \Delta$ of a $Cut_{cs}$-inference have separated derivations in ${\bf G3c[mic]}^{\cal R} +Cut_{cs}$, one of which is free of logical inferences, then its conclusion $\Gamma\seq \Delta $ has a separated derivation in the same system.
\end{lemma}

{\bf Proof}  If both derivations are free of logical rules then it suffices to apply a $Cut_{cs}$-inference to their endsequents. Otherwise we distinguish two cases.

 Case 1 ${\cal D}$ is  a derivation without logical inferences of $\Gamma \seq \Delta, A$ and ${\cal E}$ is  a separated derivation containing logical inferences  of $A, \Gamma\seq \Delta$, so that it ends with a logical inference.
 
 We have to find a separated derivation of $\Gamma\seq \Delta$.  By Lemma \ref{atomic} there is an atomic subsequent $\Gamma^\circ \seq \Delta^\circ$ of  $\Gamma\seq \Delta, A$ such that  $\Gamma^\circ \seq \Delta^\circ$ has  a derivation ${\cal D}^\circ$ without  logical inferences. If $A$ does not occur in  $\Delta^\circ$, then a separated derivation,   actually a derivation without logical inferences, of $\Gamma \seq \Delta$ can be obtained directly by weakening the conclusion of ${\cal D}^\circ$. On the other hand if $A$ occurs in $\Delta^\circ$, then $A$ is atomic so that  the principal formula of the last logical inference of ${\cal E}$  is different from $A$. We can then apply the induction hypothesis to ${\cal D}^\circ$, appropriately weakened,  and to the immediate subderivation(s) of ${\cal E}$ and then the same last logical inference of ${\cal E}$. For example if ${\cal E}$ has the form:
 \[
 \ba{c}
 {\cal E}_0\\
 A, E, F, \Gamma'\seq \Delta\\
 \cline{1-1}
 A,E\et F, \Gamma'\seq \Delta
 \ea
 \]
by Lemma \ref{atomic} we have a derivation ${\cal D}'$ without logical inferences of $E,F,\Gamma'\seq \Delta, A$.
 By induction hypothesis applied to ${\cal D}'$ and ${\cal E}_0$ we obtain a separated derivation of $E, F, \Gamma'\seq \Delta$, from which the desired separated derivation of $E\et F, \Gamma'\seq \Delta$ is obtained by means of the last $L\et$-inference of ${\cal E}$.
 
Case 2.  ${\cal D}$ is  a separated derivation containing logical inferences, so that it ends with a logical inference,  of $\Gamma \seq \Delta, A$ and ${\cal E}$ is  a  derivation without  logical inferences  of $A, \Gamma\seq \Delta$. By Lemma \ref{atomic} there is an atomic subsequent of $\Gamma^\circ\seq \Delta^\circ$ of $A,\Gamma\seq \Delta$ with a derivation ${\cal E}^\circ$ without logical inferences. If $A$ does not occur in $\Gamma^\circ$, then $\Gamma\seq \Delta$ can be derived without logical inferences by weakening the conclusion of ${\cal E}^\circ$. Otherwise $A$ is atomic, so that it is not the principal formula
  of the last inference of ${\cal D}$, and we can apply the induction hypothesis to the immediate subderivation(s) of ${\cal D}$ to  conclude as in the previous case. $\Box$

  \begin{proposition} \label{sepcutcs}

  If the  premisses $\Gamma \seq \Delta, A$  and $A, \Gamma\seq \Delta$ of a $Cut_{cs}$-inference have separated derivation in 
  ${\bf G3c[mic]}^{\cal R} + Cut_{cs}$, then its conclusion has a separated derivation in the same system.
   \end{proposition}
  
  {\bf Proof} Let ${\cal D}$ and ${\cal E}$ be  separated derivations of $\Gamma \seq \Delta, A$  and $A, \Gamma\seq \Delta$  respectively.
  We have to find a separated derivation of $\Gamma\seq \Delta$.
  
  By the previous Lemma we can assume that both ${\cal D}$ and ${\cal E}$ end with a logical rule, and  proceed by
  a principal induction  on the height (of the formation tree) of $A$ and a secondary induction  on $h({\cal D})+h({\cal E})$.

{\bf  Classical case}

 Case 1 $A$ is not principal  in (the last inference of) ${\cal D}$.
 If ${\cal D}$ reduces to an instance of $L\bot$, the same holds for $\Gamma\seq \Delta$.
 
 Case 1. $L\et$. $~ {\cal D}$ is of the form 
 
 \[
 \ba{c}
 {\cal D}_0\\
 E,F,\Gamma'\seq \Delta, A\\
 \cline{1-1}
 E\et F,\Gamma'\seq \Delta, A
 \ea
 \]so that the endsequent of ${\cal E}$ has the form $A,E\et F, \Gamma'\seq \Delta$.
 By Proposition \ref{invert} there is a separated derivation ${\cal E}'$  of $A, E, F, \Gamma'\seq\Delta$ such that $h({\cal E}')\leq h({\cal E})$.
 By the (secondary) induction hypothesis applied to ${\cal D}_0$ and ${\cal E}'$ there is a separated derivation of $E,F,\Gamma'\seq \Delta$, from which the required separated derivation of $\Gamma\seq \Delta$ can be obtained by means of the last $L\et$- inference of ${\cal E}$. In the following we will express the argument as follows: 
\[
\ba{ccc}
 E,F,\Gamma' \seq \Delta, A&~~~&\\
 \cline{1-1}
 E\et F,\Gamma'\seq \Delta, A&&A,E\et F, \Gamma'\seq \Delta\\
 \cline{1-3}
 \multicolumn{3}{c}{E\et F, \Gamma'\seq \Delta}
 \ea
 \]
 
 is transformed into:
 \[
 \ba{cccl}
&~~~&A, E\et F, \Gamma'\seq \Delta&\vbox to 0pt{\hbox{inv}}\\
 \cline{3-3}
 E, F, \Gamma'\seq\Delta, A&&A, E, F, \Gamma'\seq \Delta&\vbox to 0pt{\hbox{ind}}\\
 \cline{1-3}
 \multicolumn{3}{c}{E,F,\Gamma'\seq \Delta}\\
  \multicolumn{3}{c}{\overline{E\et F,\Gamma'\seq \Delta}}\
 \ea
 \]
In this case the  principal formula of the last logical inference of ${\cal D}$ occurs in the endsequent of ${\cal E}$ where it can be inverted, producing a separated derivation of a sequent identical to the  premiss of the last inference of ${\cal D}$, except that the cut formula $A$ is shifted from the succedent to the antecedent.
We can then apply the secondary induction hypothesis to produce  a sequent to which the last logical inference of ${\cal D}$ can be applied, yielding the required separated derivation.
The same kind of  argument applies to all the remaining cases. For example:

 Case 1. $L\imp$
 
 \[
 \ba{ccccc}
 \Gamma'\seq \Delta, E, A&~~&F,\Gamma'\seq \Delta, A&~~~&\\
 \cline{1-3}
 \multicolumn{3}{c}{E\imp F, \Gamma'\seq \Delta, A}&&A,E\imp F, \Gamma'\seq \Delta\\
 
  \multicolumn{5}{c}{\overline{~~~~~~~~~~~~~~~~~~~ E\imp F, \Gamma'\seq \Delta~~~~~~~~~~~~~~~~~~~~~~~~~~}}
   \ea
 \]
 is transformed into:
 
 \[
 \ba{ccclcccl}
 &~~&A,E\imp F, \Gamma'\seq \Delta&\vbox to 0pt{\hbox{inv}}~~~~&&~~&A,E\imp F, \Gamma'\seq \Delta&\vbox to 0pt{\hbox{inv}}\\
 \cline{3-3}\cline{7-7}
 \Gamma'\seq \Delta, E, A&&A,\Gamma'\seq \Delta, E&\vbox to 0pt{\hbox{ind}}~~~~~~~& F,\Gamma'\seq \Delta, A&&A,F,\Gamma'\seq \Delta&\vbox to 0pt{\hbox{ind}}\\
 \cline{1-3}\cline{5-7}
 \multicolumn{3}{c}{\Gamma'\seq \Delta, E}&&\multicolumn{3}{c}{F,\Gamma' \seq \Delta}\\
 \multicolumn{7}{c}{\overline{~~~~~~~~~~~~~~~~~~~~~E \imp F, \Gamma' \seq \Delta~~~~~~~~~~~~~~~~~~~~~~~~~~}}
 \ea
 \]

 Case1 $R\imp$
 \[
 \ba{ccc}
 \Gamma, E\seq \Delta', F, A&~~&\\
 \cline{1-1}
 \Gamma\seq \Delta', E\imp F, A&&A,\Gamma\seq \Delta', E\imp F\\
 \cline{1-3}
  \multicolumn{3}{c}{\Gamma\seq \Delta', E\imp F}
 \ea
 \]
 is transformed into:
 
 \[
\ba{cccl}
&&A,\Gamma\seq \Delta', E\imp  F&\vbox to 0pt{\hbox{inv}}\\
\cline{3-3}
\Gamma,E\seq \Delta', F, A&&A,\Gamma, E  \seq \Delta',  F&\vbox to 0pt{\hbox{ind}}\\
\cline{1-3}
 \multicolumn{3}{c}{\Gamma, E\seq \Delta',   F}\\
 \multicolumn{3}{c}{\overline{\Gamma\seq \Delta', E\imp F}}
\ea
\]

 Case 1 $R\forall$
 
 \[
 \ba{ccc}
 \Gamma\seq \Delta', E[x/a], A&~~&\\
 \cline{1-1}
 \Gamma \seq \Delta', \forall x E, A&&A,\Gamma\seq \Delta', \forall x E\\
 \cline{1-3}
 \multicolumn{3}{c}{ \Gamma\seq \Delta', \forall x E}
 \ea
 \]
 
 is transformed into:
 
 \[
 \ba{clcl}
&~~&A,\Gamma\seq \Delta', \forall x E&\vbox to 0pt{\hbox{inv}}\\
\cline{3-3}
 \Gamma\seq \Delta', E[x/a], A &&A,\Gamma\seq \Delta', E[x/a]& \vbox to 0pt{\hbox{ind}}\\
 \cline{1-3}
 \multicolumn{3}{c}{ \Gamma\seq \Delta', E[x/a]}\\
 \multicolumn{3}{c}{\overline{ \Gamma\seq \Delta', \forall x E}}
 \ea
 \]

 Case 2.  $A$ is not principal  in (the last inference of) ${\cal E}$.
 If ${\cal E}$ reduces to an instance of $L\bot$, the same holds for $\Gamma\seq \Delta$.
 All the other cases are treated dually to Case 1. For example:
 
 Case 2. $R\et$
 
 \[
 \ba{ccc}
 &~~&A,E,F,\Gamma'\seq \Delta\\
 \cline{3-3}
 E\et F, \Gamma'\seq \Delta, A&&A,E\et F, \Gamma'\seq \Delta\\
 \cline{1-3}
 \multicolumn{3}{c}{E\et F, \Gamma'\seq \Delta}
  \ea
 \]
  is transformed into:
 \[
 \ba{clcl}
 E\et F, \Gamma'\seq \Delta, A&\vbox to 0pt{\hbox{inv}}~~&&\\
 \cline{1-1}
 E,F,\Gamma'\seq \Delta, A&&A,E,F,\Gamma'\seq \Delta&\vbox to 0pt{\hbox{ind}}\\
 \cline{1-3}
 \multicolumn{3}{c}{E,F, \Gamma'\seq \Delta}&\\
 \multicolumn{3}{c}{\overline{E\et F, \Gamma'\seq \Delta}}
 \ea
 \]
In   the above cases all  the applications of $\hbox{ind}$ refer to  the secondary induction hypothesis.
 
 \

 Case 3. $A$ is principal  in (the last inferences of) both $ {\cal D}$  and  ${\cal E}$.
 
 \

 Case 3. $\et$:
 
 \[
 \ba{ccccc}
 \Gamma\seq \Delta, B&~~~&\Gamma\seq \Delta, C&~~~&B, C,\Gamma\seq \Delta\\
 \cline{1-3}\cline{5-5}
 \multicolumn{3}{c}{\Gamma \seq \Delta, B\et C}&&B\et C,\Gamma\seq \Delta\\
 \cline{1-5}
  \multicolumn{5}{c}{\Gamma \seq \Delta}
 
 \ea
 \]
is transformed into:

\[
\ba{ccclcl}
&~~&\Gamma\seq \Delta, B&\vbox to 0pt{\hbox{w}}~~~&&\\
\cline{3-3}
&&C,\Gamma\seq \Delta, B&&B,C,\Gamma\seq \Delta&\vbox to 0pt{\hbox{ind}}\\
\cline{3-5}
\Gamma\seq \Delta, C&&\multicolumn{3}{c}{C,\Gamma\seq \Delta}&\vbox to 0pt{\hbox{ind}}\\
\cline{1-5}
 \multicolumn{5}{c}{\Gamma \seq \Delta}
\ea
\]

In this case  the second application of $\hbox{ind}$  refers necessarily to the principal induction hypothesis, which is possible since $h(C)< h(B\et C)$,
independently of the height of the separated derivation of the second premiss $C,\Gamma\seq \Delta$,  previously obtained by the (secondary suffices) induction hypothesis. Of the remaining cases we deal with the  Case 3 $\forall$  in which $A$ is a universal formula, leaving the others to the reader.

\[
\ba{ccc}
\Gamma\seq \Delta, B[x/a]&~~~&\forall x B, B[x/t],\Gamma\seq \Delta\\
\cline{1-1}\cline{3-3}
\Gamma\seq \Delta,\forall x B&&\forall x B,\Gamma\seq \Delta\\
\cline{1-3}
\multicolumn{3}{c}{\Gamma\seq \Delta}
\ea
\]
 is transformed into:

 \[
 \ba{clclcl}
 &&\Gamma\seq \Delta, B[x/a]&&&\\
 \cline{3-3}
 &&\Gamma\seq \Delta,\forall x B&\vbox to 0pt{\hbox{w}}~~&&\\
 \cline{3-3}
 \Gamma\seq \Delta, B[x/a]&\vbox to 0pt{\hbox{Sub$[a/t]$}}~~&B[x/t],\Gamma\seq \Delta, \forall x B&&\forall x B,B[x/t],\Gamma\seq \Delta&\vbox to 0pt{\hbox{ind}}\\
 \cline{1-1}\cline{3-5}
 \Gamma\seq \Delta, B[x/t]&&\multicolumn{3}{c}{B[x/t], \Gamma\seq \Delta}&\vbox to 0pt{\hbox{ind}}\\
 \cline{1-5}
 \multicolumn{5}{c}{\Gamma\seq \Delta}
 \ea
 \]where $\hbox{Sub$[a/t]$}$ yields the result of replacing $a$ by $t$ throughout the separated derivation of $\Gamma\seq \Delta, B[x/a]$.
For the result of  such a replacement to be a derivation it is required that the parameters used as proper in the $L\exists$ and $R\forall$-inferences of  the given derivation be renamed, if necessary, so as not to occur in $t$.

{\bf Intuitionistic case}

Case 1. $A$ is not principal in ${\cal D}$.
We only need to replace Case1. $L\imp$, $R\imp$ and   $R\forall$ with the following:

Case 1. $L^i\imp$
\[
\ba{ccccc}
E\imp F, \Gamma'\seq \Delta, E,  A&~~&F,\Gamma'\seq \Delta, A&~~&\\
\cline{1-3}
\multicolumn{3}{c}{E\imp F, \Gamma'\seq \Delta, A}&&A,E\imp F,\Gamma'\seq \Delta\\
\cline{2-5}
&\multicolumn{4}{c}{E\imp F, \Gamma'\seq \Delta}
\ea
\]
is transformed into:

\[
\ba{ccclcccl}
&&A,E\imp F,\Gamma'\seq \Delta&\vbox to 0pt{\hbox{inv}}~~&&&A,E\imp F,\Gamma'\seq \Delta&\vbox to 0pt{\hbox{inv}}\\
\cline{3-3}\cline{7-7}
E\imp F, \Gamma'\seq \Delta, E,  A&~~&A,E\imp F,\Gamma'\seq \Delta, E&\vbox to 0pt{\hbox{ind}}&F,\Gamma'\seq \Delta, A&&A, F, \Gamma'\seq \Delta&\vbox to 0pt{\hbox{ind}}\\
\cline{1-3}\cline{5-7}
\multicolumn{3}{c}{E\imp F,\Gamma'\seq \Delta, E}&& \multicolumn{3}{c}{F,\Gamma'\seq \Delta}\\
\multicolumn{7}{c}{\overline{~~~~~~~~~~~~~~~~~~~~~~~~~~~~~~~ E\imp F, \Gamma'\seq \Delta ~~~~~~~~~~~~~~~~~~~~~~~~~~~~~~~~~~~~~~}}
\ea
\]

Case 1. $R^i\imp$

\[
\ba{ccc}
E,\Gamma\seq F&~~~&\\
\cline{1-1}
\Gamma\seq \Delta', E\imp F, A&&A,\Gamma\seq \Delta',E\imp F\\
\cline{1-3}
\multicolumn{3}{c}{\Gamma\seq \Delta',  E\imp F}
\ea
\]

is transformed into:
\[
\ba{c}
E,\Gamma\seq F\\
\cline{1-1}
\Gamma\seq \Delta', E\imp F
\ea
\]

Case 1. $R^i\forall$.

\[
\ba{ccc}
\Gamma\seq E[x/a]&~~&\\
\cline{1-1}
\Gamma\seq \Delta', \forall x E, A&&A,\Gamma\seq \Delta', \forall x E\\
\cline{1-3}
\multicolumn{3}{c}{\Gamma\seq \Delta', \forall x E}

\ea
\]

is transformed into:

\[
\ba{c}
\Gamma\seq E[x/a]\\
\cline{1-1}
\Gamma\seq \Delta', \forall x E
\ea
\]

Case 2  Case 1. does not occur, so that $A$ is principal in the last inference of ${\cal D}$,  and $A$   is not principal in ${\cal E}$ and
If the rule  of the last inference of ${\cal E}$ is invertible in ${\cal D}$ , i.e. if it is not a $R^i\imp$ or $R^i\forall$-inference then we proceed as in Case 2 of the classical case.
If ${\cal E}$ ends with a $R^i\imp$ or $R^i\forall$-inference,  we distinguish cases according to the form of $A$.

\

Case 2. $R^i\imp$, $\et$

\[
\ba{ccccc}
\Gamma\seq \Delta', E\imp F, B&~~&\Gamma\seq \Delta', E\imp F, C&~~&B\et C, E,\Gamma\seq F\\
\cline{1-3} \cline{5-5}
\multicolumn{3}{c}{\Gamma\seq \Delta', E\imp F, B\et C}&&B\et C, \Gamma\seq \Delta', E\imp F\\
\multicolumn{5}{c}{\overline{~~~~~~~~~~~~~~~~~~~~~~~~~~~\Gamma\seq \Delta', E\imp F~~~~~~~~~~~~~~~~~~~~~~~~~~~}}
\ea
\]
is transformed into:

\[
\ba{ccclll}
&&&&B\et C, E,\Gamma\seq F&\\
\cline{5-5}
&~~&\Gamma\seq\Delta', E\imp F, B&\vbox to 0pt{\hbox{w}}~~&B\et C, \Gamma\seq \Delta', E\imp F&\vbox to 0pt{\hbox{inv}}\\
\cline{3-3}\cline{5-5}
&&C,\Gamma\seq\Delta', E\imp F, B&&B,  C, \Gamma\seq \Delta', E\imp F& \vbox to 0pt{\hbox{ind}}\\
\cline{3-5}
\Gamma\seq \Delta',E\imp F, C&\multicolumn{4}{c}{~~~~~~~~~~~~C,\Gamma\seq \Delta', E\imp F~~}& \vbox to 0pt{\hbox{ind}}\\
\cline{1-5}
\multicolumn{4}{c}{\Gamma\seq \Delta', E\imp F}
\ea
\]

Case 2. $R^i\imp$, $\vel$

\[
\ba{ccc}
\Gamma\seq \Delta', E\imp F, B, C&~~~&B\vel C, E,\Gamma \seq F\\
\cline{1-1}\cline{3-3}
\Gamma\seq \Delta', E\imp F, B \vel C&&B\vel C, \Gamma \seq \Delta',  E\imp F\\
\cline{1-3}
\multicolumn{3}{c}{\Gamma\seq \Delta', E\imp F}

\ea
\]
is transformed into:

\[
\ba{ccclcl}
&&B\vel C, E, \Gamma\seq F&&&\\
\cline{3-3}
&~~&B\vel C, \Gamma\seq \Delta', E\imp F&\vbox to 0pt{\hbox{inv}}~~&&\\
\cline{3-3}
&& C, \Gamma\seq \Delta', E\imp F&\vbox to 0pt{\hbox{w}}&B\vel C, E,\Gamma \seq F&\\
\cline{3-3} \cline{5-5}
\Gamma\seq \Delta', E\imp F, B, C&&C, \Gamma\seq \Delta', E\imp F, B&\vbox to 0pt{\hbox{ind}}&B\vel C, \Gamma\seq \Delta', E\imp F&\vbox to 0pt{\hbox{inv}}\\
\cline{1-3}\cline{5-5}
\multicolumn{3}{c}{\Gamma\seq \Delta', E\imp F, B}&&B , \Gamma\seq \Delta', E\imp F&\vbox to 0pt{\hbox{ind}}\\

\multicolumn{6}{c}{\overline{~~~~~~~~~~~~~~~~~~~~~~\Gamma\seq \Delta', E\imp F~~~~~~~~~~~~~~~~~~~~~~~~~~~~~~~~~~~}}

\ea
\]

Case 2 $R^i\imp$, $\imp$
\ 

\[
\ba{ccc}
B,\Gamma\seq C&~~~&B\imp C, E,\Gamma\seq F\\
\cline{1-1}\cline{3-3}
\Gamma\seq  \Delta', E\imp F,B\imp C&&B\imp C,\Gamma\seq \Delta',E\imp F\\
\cline{1-3}
\multicolumn{3}{c}{\Gamma\seq \Delta', E\imp F}
\ea
\]
is transformed into:

\[
\ba{clcl}
B,\Gamma\seq C&~~~&&\\
\cline{1-1}
\Gamma\seq F,  B\imp C&\vbox to 0pt{\hbox{w}}&&\\
\cline{1-1}
E,\Gamma\seq F, B\imp C&&B\imp C, E, \Gamma\seq F&\vbox to 0pt{\hbox{ind}}\\
\cline{1-3}
\multicolumn{3}{c}{E, \Gamma\seq  F}\\
\multicolumn{3}{c}{\overline{~~ \Gamma\seq  \Delta',E\imp  F}}

\ea
\]

Case 2 $R^i\imp$, $\forall$

\[
\ba{ccc}
\Gamma\seq B[x/a]&~~~&\forall x B, E, \Gamma\seq F\\
\cline{1-1}\cline{3-3}
\Gamma\seq \Delta', E\imp F, \forall x B&&\forall x B, \Gamma\seq \Delta', E\imp F\\
\cline{1-3}
\multicolumn{3}{c}{~~ \Gamma\seq  \Delta',E\imp  F}
\ea
\]
is transformed into:
\[
\ba{clcl}
\Gamma\seq B[x/a]&&&\\
\cline{1-1}
\Gamma\seq F, \forall x B &\vbox to 0pt{\hbox{w}}~~&&\\
\cline{1-1}
E,\Gamma\seq F, \forall x B&&\forall x B, E, \Gamma\seq F&\vbox to 0pt{\hbox{ind}}\\
\cline{1-3}
\multicolumn{3}{c}{E,  \Gamma\seq    F}\\
\multicolumn{3}{c}{\overline{~~ \Gamma\seq  \Delta',E\imp  F}}
\ea
\]

Case 2 $R^i\imp$, $\exists$

\ 

\[
\ba{ccc}
\Gamma\seq \Delta', E\imp F, \exists xB, B[x/t]&~~&\exists x B, E,\Gamma\seq  F\\
\cline{1-1}\cline{3-3}
\Gamma\seq \Delta', E\imp F, \exists x B&& \exists x B, \Gamma\seq \Delta', E\imp F\\
\cline{1-3}
\multicolumn{3}{c}{\Gamma\seq  \Delta', E\imp  F}
\ea
\]
is transformed into:

\[
\ba{ccclcl}
&&&&\exists x B, E,\Gamma\seq F&\\
\cline{5-5}
&~~&\exists x B, E,\Gamma\seq  F&~~~~&  \exists x B, \Gamma\seq \Delta', E\imp F&\vbox to 0pt{\hbox{inv}}\\
\cline{3-3}\cline{5-5}
\Gamma\seq \Delta', E\imp F, \exists x B, B[x/t]&& \exists x B, \Gamma\seq \Delta',E\imp  F, B[x/t]&\vbox to 0pt{\hbox{ind}}&B[x/a], \Gamma\seq \Delta', E\imp F&\vbox to 0pt{\hbox{Sub$[a/t]$}}\\
\cline{1-3}\cline{5-5}
\multicolumn{3}{c}{\Gamma\seq \Delta', E\imp F, B[x/t] } && B[x/t], \Gamma\seq \Delta', E\imp F&\vbox to 0pt{\hbox{ind}}\\
\multicolumn{6}{c}{\overline{~~~~~~~~~~~~~~~~~~~~~~~~~~\Gamma\seq  \Delta',E\imp  F~~~~~~~~~~~~~~~~~~~~~~~~~~~~~~~~~~~~~~~~~~~~~~~}}

\ea
\]

Case 2. $R^i\forall$.  Similar to Case2. $R^i\imp$. For example:
 
 (Case 2 $R^i\forall$, $\imp$)
 
 \[
 \ba{ccc}
 B,\Gamma\seq C&~~~&B\imp C,\Gamma\seq E[x/a]\\
 \cline{1-1}\cline{3-3}
 \Gamma\seq \Delta', \forall x E, B\imp C&&B\imp C, \Gamma\seq \Delta', \forall x E\\
 \cline{1-3}
 \multicolumn{3}{c}{\Gamma\seq  \Delta',\forall x E}
 
 \ea
 \]
is transformed into:
\[
\ba{cccl}
 B,\Gamma\seq C&~~~&&\\
 \cline{1-1}
 \Gamma\seq E[x/a], B\imp C&&B\imp C,\Gamma\seq E[x/a]&\vbox to 0pt{\hbox{ind}}\\
 \cline{1-3}
 
 \multicolumn{3}{c}{\Gamma\seq   E[x/a]}\\
  \multicolumn{3}{c}{\overline{\Gamma\seq   \Delta', \forall x E}}
\ea
\]

and  (Case 2  $R^i\forall$, $\forall$)

\[
\ba{ccc}
\Gamma\seq B[x/b]&~~~&\forall x B,\Gamma\seq E[x/a]\\
\cline{1-1}\cline{3-3}
\Gamma\seq \Delta',\forall x E,\forall x B&&\forall x B,\Gamma\seq \Delta',\forall x E\\
\cline{1-3}
 \multicolumn{3}{c}{\Gamma\seq  \Delta', \forall x  E}

\ea
\]
is transformed into:

\[
\ba{cccl}
\Gamma\seq B[x/b] &~~&&\\
\cline{1-1}
\Gamma\seq E[x/a], \forall x B&&\forall x B, \Gamma\seq E[x/a]&\vbox to 0pt{\hbox{ind}}\\
\cline{1-3}
\multicolumn{3}{c}{\Gamma\seq   E[x/a]}\\
  \multicolumn{3}{c}{\overline{\Gamma\seq   \Delta', \forall x E}}
\ea
\]

Case 3. $A$ is principal in both  ${\cal D}$ and ${\cal E}$. The only difference  with respect to the classical case concern $\imp$ and $\forall$.

Case 3i $\imp$

\[
\ba{ccccc}
B,\Gamma\seq C&~~&B\imp C, \Gamma\seq \Delta, B&~~&C,\Gamma\seq \Delta\\
\cline{1-1}\cline{3-5}
\Gamma\seq \Delta, B\imp C&&\multicolumn{3}{c}{B\imp C,\Gamma \seq \Delta}\\ 
\cline{1-5}
\multicolumn{5}{c}{\Gamma\seq  \Delta}
\ea
\]
is transformed into:

\[
\ba{ccclcccl}
B,\Gamma\seq C&~~&&&B,\Gamma\seq C&\vbox to 0pt{\hbox{w}} & C,\Gamma\seq \Delta&\vbox to 0pt{\hbox{w}}\\
\cline{1-1}  \cline{5-5} \cline{7-7}
\Gamma\seq \Delta, B,  B\imp C&~~&B\imp C, \Gamma\seq \Delta, B&  \vbox to 0pt{\hbox{ind}}~~~~& B,\Gamma\seq \Delta, C&~~&C,B,\Gamma\seq \Delta&\vbox to 0pt{\hbox{ind}}\\
\cline{1-3}\cline{5-7}
\multicolumn{3}{c}{\Gamma\seq \Delta, B}&&\multicolumn{3}{c}{B, \Gamma\seq \Delta}&\vbox to 0pt{\hbox{ind}}\\
\cline{1-7}
\multicolumn{7}{c}{\Gamma\seq \Delta}
\ea
\]

Case 3i $\forall$

\[
\ba{ccc}
\Gamma\seq B[x/a]&~~~&\forall x B, B[x/t],\Gamma\seq \Delta\\
\cline{1-1}\cline{3-3}
\Gamma\seq \Delta,\forall x B&&\forall x B,\Gamma\seq \Delta\\
\cline{1-3}
\multicolumn{3}{c}{\Gamma\seq \Delta}

\ea
\]
is transformed into:

\[
\ba{clclclcl}
&&\Gamma\seq B[x/a]&&&\\
\cline{3-3}
\Gamma\seq B[x/a]&\vbox to 0pt{\hbox{Sub$[a/t]$}}~~~~~&\Gamma\seq \Delta, \forall x B&\vbox to 0pt{\hbox{w}}&&\\
\cline{1-1}\cline{3-3}
\Gamma\seq B[x/t]&\vbox to 0pt{\hbox{w}}~~~&B[x/t], \Gamma\seq \Delta, \forall x B&~~&\forall x B,B[x/t],\Gamma\seq \Delta&\vbox to 0pt{\hbox{ind}}\\
\cline{1-1}\cline{3-5}
\Gamma\seq \Delta, B[x/t]&&\multicolumn{3}{c}{B[x/t], \Gamma\seq \Delta}&\vbox to 0pt{\hbox{ind}}\\
\cline{1-5}
\multicolumn{5}{c}{\Gamma\seq \Delta}

\ea
\]
$\Box$

From  Proposition \ref{sepequality} and \ref{sepcutcs} by a straightforward induction argument we have the following separation  property
for ${\bf G3[mic]}^{\cal R} + Cut_{cs} $.

\begin{proposition}\label{sep}

Every derivation in ${\bf G3[mic]}^{\cal R} + Cut_{cs} $ can be transformed into a separated derivation of its endsequent.

\end{proposition}

\begin{theorem}\label{teorsep}
If the structural rules  are admissible in ${\cal R}$,  then they are admissible in ${\bf G3[mic]}^{\cal R}$ as well.
\end{theorem}

{\bf Proof}
 Let   ${\cal D}$ be a derivation in ${\bf G3[mic]}^{\cal R} + RLW +RLC+ Cut $.  We have to show that the applications of the $RLW$, $RLC$ and $Cut$ can be eliminated from ${\cal D}$.
  ${\cal D}$   can be transformed into a derivation ${\cal D}'$
in  ${\bf G3[mic]}^{\cal R} + Cut_{cs} $  of the same endsequent.
For,   
the application of the $Cut$-rule can be replaced by applications of the weakenings and the $Cut_{cs}$ rule.

The applications of the contraction rules can be replaced by derivations from its premiss and initial sequents using  the $Cut_{cs}$-rule. More precisely,  as far as left contraction is concerned, the subderivations of ${\cal D}$ of the form
\[
\ba{c}
{\cal E}\\
\underline{F, F, \Gamma\seq \Delta}\\
F, \Gamma\seq \Delta
\ea
\]
can be replaced  by:

\[
\ba{ccc}
{\cal I} &~~~&{\cal E}\\
F,\Gamma\seq F&&F,F,\Gamma\seq \Delta\\
\cline{1-3}
\multicolumn{3}{c}{F,\Gamma\seq \Delta}
\ea
\]where, in case $F$ is not atomic, ${\cal I}$ is a derivation in ${\bf G3m}$ or in ${\bf G3i}$. 
Similarly for right contraction.

Finally the applications of the weakening rules  can be eliminated by their (heigh-preserving) admissibiity in  all the systems considered.

Thus from ${\cal D}$ we obtain a derivation ${\cal D}'$ in ${\bf G3[mic]}^{\cal R} +Cut_{cs}$, that 
by Proposition \ref{sep},  can be transformed into 
a separated derivation ${\cal D}''$ of the endsequent  of ${\cal D}$. Thus,  to obtain the desired derivation in ${\bf G3}[mic]^{\cal R}$ of the endsequent of ${\cal D}$,  it suffices to eliminate the applications of $Cut_{cs}$ in the initial subderivations of ${\cal D}''$ belonging to ${\cal R} +Cut_{cs}$, which  is possible if the  contraction and the cut rule are admissible in ${\cal R}$.
$\Box$

\subsection{Admissibility of the Structural Rules}

\begin{corollary}
The structural rules are admissible in ${\bf G3[mic]}$
\end{corollary}

{\bf Proof} By Theorem \ref{teorsep},  with ${\cal R}=\es$,  
it suffices to note that the sequents that  can be derived  from
initial sequents by means of the structural rules are themselves initial sequents. $\Box$

\ 

Following \cite{N16} we let 
${\bf G3[mic]}^=$ be   ${\bf G3[mic]}^{\cal R}$  for ${\cal R}=\{\hbox{Ref}, \hbox{Repl}\}$.

\begin{corollary}
The structural rules are admissible in ${\bf G3[mic]}^=$
\end{corollary}
 
 {\bf Proof} 
 For the admissibility of left contraction we proceed by induction on the height of   derivations  to show that a derivation ${\cal D}$ of $A,A,\Gamma\seq \Delta$ in ${\cal R}$ can be transformed into a derivation of $A,\Gamma \seq \Delta$ in ${\cal R}$. That is immediate if ${\cal D}$ reduces to an initial sequent.  If $h({\cal D})>0$,  the conclusion is a straightforward consequence of the induction hypothesis except when
  ${\cal D}$ has the form:
  
 \[
 \ba{cl}
 \underline{s=r, s=r, E[x/r], \Gamma'\seq \Delta}&\vbox to 0pt{\hbox{Repl}}\\
 s=r, s=r, \Gamma' \seq \Delta
 \ea
 \] and $A$ is $s=r$,  that coincides also  with $E[x/s]$.
We may assume that there is exactly one occurrence of $x$ in $E$.
Then $E$ can have the form $x=r$ or $r$ can have the form $r^\circ[x/s]$ and $E$ the form $s=r^\circ$.
In the former case the induction hypothesis applied to the immediate subderivation ${\cal D}_0$ of  ${\cal D}$,  whose endsequent is 
$s=r,s=r,r=r, \Gamma'\seq \Delta$, yields a derivation 
of $s=r, r=r, \Gamma'\seq \Delta$, from which we obtain $s=r, \Gamma'\seq \Delta$ by an application of $\hbox{Ref}$.
In the latter case ${\cal D}$ has the form:

\[
\ba{cl}
{\cal D}_0\\
\underline{s=r^\circ[x/s], ~s=r^\circ[x/s],~ s=r^\circ[x/r^\circ[x/s]], ~\Gamma'\seq \Delta}&\vbox to 0pt{\hbox{Repl}}\\
s=r^\circ[x/s],~ s=r^\circ[x/s], ~\Gamma'\seq \Delta

\ea
\]
and  can be transformed into:

 \[
\ba{cl}
{\cal D}_0&\\
\underline{s=r^\circ[x/s], ~s=r^\circ[x/s],~ s=r^\circ[x/r^\circ[x/s]], ~\Gamma'\seq \Delta}&\vbox to 0pt{\hbox{Repl}}\\
s=r^\circ[x/s],~ s=r^\circ[x/s], ~\Gamma'\seq \Delta&\vbox to 0pt{\hbox{w}}\\
\overline{\underline{s=r^\circ[x/s],~s=s,~ s=r^\circ[x/s], ~\Gamma'\seq \Delta}}&\vbox to 0pt{\hbox{Repl}}\\
\underline{s=r^\circ[x/s],~s=s, ~\Gamma'\seq \Delta}&\vbox to 0pt{\hbox{Ref}}\\
s=r^\circ[x/s], ~\Gamma'\seq \Delta
\ea
\]
 Since the rules of ${\cal R}$ do not modify the succedent of their premiss, the admissibility of the cut rule follows  by a straightforward induction on the height of the derivation of its first   premiss. $\Box$
 
%\subsection{Aknowledgment} We wish to thank Sara Negri for helpful comments and suggestions.

\

\end{document}